\documentclass[12pt]{amsart}
\input{epsf}
\usepackage{amssymb}
\usepackage{amsfonts}
\usepackage{amsmath}
\usepackage{amscd}

\begin{document}
\title[IH of Toric Varieties and a Conjecture of Kalai]{Intersection Homology of Toric Varieties and a Conjecture of Kalai}
\author{Tom Braden}
\address{IAS, Princeton}
\email{braden@math.ias.edu}
\thanks{The first author was supported by NSF grant DMS 9304580}
\author{Robert MacPherson}
\address{IAS, Princeton}
\maketitle

\renewcommand{\baselinestretch}{1,0}

\theoremstyle{plain}
\newtheorem{thm}{Theorem}
\newtheorem{lemma}[thm]{Lemma}
\newtheorem{prop}[thm]{Proposition}
\newtheorem{cor}[thm]{Corollary}
\theoremstyle{definition}
\newtheorem*{conj}{Conjecture}
\newtheorem*{defn}{Definition}
\newtheorem*{nota}{Notation}

\hyphenation{Mac-Pher-son}
\hyphenation{canon-ical}

\newcommand{\ra}{\rightarrow}
\newcommand{\la}{\leftarrow}
\newcommand{\lra}{{\longrightarrow}}
\newcommand{\lla}{{\longleftarrow}}
\newcommand{\wt}{\widetilde}

\def\rank{\mathop{\rm rank}\nolimits}
\def\ker{\mathop{\rm Ker}\nolimits}
\def\coker{\mathop{\rm Coker}\nolimits}
\def\im{\mathop{\rm Im}\nolimits}
\def\Id{{\rm Id}}
\def\id{{\rm id}}
\def\lcm{\mathop{\rm lcm}\limits}
\def\Hom{\mathop{\rm Hom}}

\newcommand{\C}{{\mathbb C}}
\newcommand{\Z}{{\mathbb Z}}
\newcommand{\R}{{\mathbb R}}
\newcommand{\Rge}{\R_{\ge 0}}
\newcommand{\PP}{{\mathbb P}}
\newcommand{\Q}{{\mathbb Q}}
\newcommand{\F}{{\mathbb F}}
\newcommand{\HH}{{\mathbb H}}
\newcommand{\cal}{\mathcal}
\newcommand{\cS}{{\cal S}}
\newcommand{\cF}{{\cal F}}
\newcommand{\cE}{{\cal E}}
\newcommand{\cL}{{\cal L}}
\newcommand{\cP}{{\cal P}}
\newcommand{\cQ}{{\cal Q}}
\newcommand{\cR}{{\cal R}}
\newcommand{\cD}{{\cal D}}
\newcommand{\gb}{{\bar g}}
\newcommand{\IC}{{\mathbf {IC^{\textstyle \cdot}}}}
\newcommand{\Aff}{\mathop{{\mathcal A}{\it ff}}}
\newcommand{\mb}{\mathbf}

Suppose that a $d$-dimensional convex polytope $P\subset \R^d$ is 
{\em rational}, i.e.\ its vertices 
are all rational points.  Then $P$ gives rise to a polynomial  
$g(P)=1+g_1(P)q+g_2(P)q^2+\cdots$ with non-negative coefficients as  
follows.  Let $X_P$ be the {\em associated toric variety} (see \S 6
-- our variety $X_P$ is $d+1$-dimensional and affine).   The  
coefficient $g_i$ is the rank of the $2i$-th intersection homology group of  
$X_P$.

The polynomial $g(P)$ turns out to depend only on the face lattice  
of $P$, (see \S 1).  It can be thought of as a measure of the complexity
of $P$;  for example, $g(P)=1$ if and only if $P$ is a simplex.  

Suppose that $F\subset P$ is a face of dimension $k$.  We construct an  
associated polyhedron $P/F$ as follows (see the figure below):  
choose an $(n-k-1)$-plane $L$  
whose intersection with $P$ is a single point $p$ of the interior of $F$.   
Let $L'$ be a small parallel displacement of $L$ that intersects the  
interior of $P$.  Then $P/F$ is the intersection of $P$ with $L'$;
it is only well-defined up to a projective transformation, but its
combinatorial type is well-defined.  Faces of $P/F$ are in one-to-one
correspondence with faces of $P$ which contain $F$.

\begin{figure}[h]
\begin{center}
\leavevmode
\hbox{%
\epsfxsize = 3.1in
\epsffile{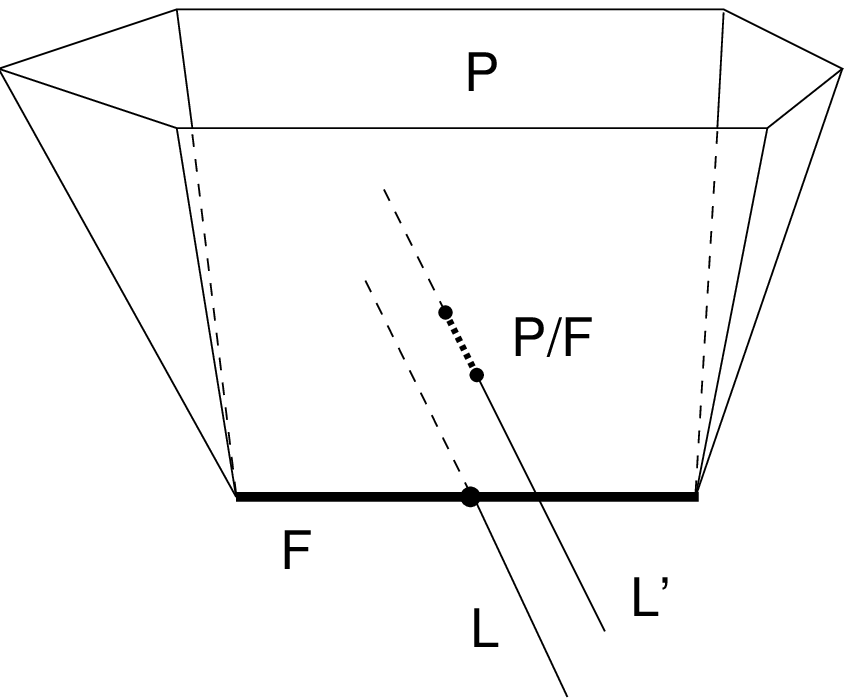}}
\end{center}
\label{fig}
\end{figure}
 
In Corollary 6, we show that
\[g(P)\geq g(F)g(P/F)\]
holds, coefficient by coefficient. 
This was conjectured by Kalai in \cite{K2}, where some of its  
applications were discussed.  The special case of the linear and 
quadratic terms was proved in \cite{K3}.  
Roughly, this inequality means that the  
complexity of $P$ is bounded from below by the complexity of the face  
$F$ and the normal complexity $g(P/F)$ to the face $F$.   
 
The principal idea is to introduce 
{\em relative $g$-polynomials $g(P,F)$} for  
any face $F$ of $P$ (\S 2).  These generalize the ordinary $g$-polynomials  
since $g(P,P)=g(P)$.  They are also combinatorially determined by the  
face lattice.  They measure the complexity of $P$ relative to the  
complexity of $F$.  For example, if $P$ is the join of $F$ with another  
polytope, then $g(P,F)=1$ (the converse, however, does not hold).
 
Our main result gives an interpretation of the coefficients 
$g_i(P,F)$ of the relative $g$-polynomials as dimensions of 
vector spaces arising  
from the topology of the toric variety $X_P$.  This shows that the  
coefficients are positive.  Kalai's conjecture is a corollary.   
 
The combinatorial definition of the relative $g$-polynomials 
$g(P,F)$ makes sense whether or not the polytope $P$ is rational.  We  
conjecture the positivity $g(P,F)\geq 0$ in general.  This would imply  
Kalai's conjecture for general polytopes.   
 
This paper is organized as follows:  The first three sections are  
entirely about the combinatorics of polyhedra.  They develop the  
properties of relative $g$-polynomials as combinatorial objects, with  
the application to Kalai's conjecture.  The last three sections concern  
algebraic geometry.  A separate guide to their contents is included in  
the introduction to \S\S 4 - 6. 
\section{$g$-numbers of polytopes}

Let $P \subset \R^d$ be a $d$-dimensional convex polytope, i.e.\ the convex
hull of a finite collection of points affinely spanning $\R^d$.
The set of faces of $P$, ordered by inclusion, forms a
poset which we will denote by $\cF(P)$.  We include the empty
face $\emptyset = \emptyset_P$ and $P$ itself as members of $\cF(P)$.
It is a graded poset, with the grading given by the dimension of
faces.  By convention we set $\dim \emptyset = -1$.  Faces of $P$ of
dimension $0$, $1$, and $d-1$ will be referred to as vertices, edges,
and facets, respectively.

Given a face $F$ of $P$, the poset $\cF(F)$ is clearly isomorphic to
the interval $[\emptyset, F] \subset \cF(P)$.  The interval $[F, P]$ is
the face poset of the polytope $P/F$ defined in the introduction.  

Given the polytope $P$, there are associated polynomials (first introduced 
in \cite{S1})
$g(P) = \sum g_i(P)q^i$ and $h(P) = \sum h_i(P)q^i$, defined recursively as follows:
\begin{itemize}
\item $g(\emptyset) = 1$
\item $h(P) = \Sigma_{\emptyset \le F < P} (q-1)^{\dim P - \dim F - 1} g(F)$, and
\item $g_0(P) = h_0(P)$, $g_i(P) = h_i(P) - h_{i-1}(P)$ for $0<i\le \dim P/2$, 
and $g_i(P) = 0$ for all other $i$.
\end{itemize}
The coefficients of these polynomials will be referred to as the
$g$-numbers and $h$-numbers of $P$, respectively.  For our purposes,
the $g$-polynomials will be of primary interest; the $h$-polynomials
will not play a role here.

These numbers clearly depend only on the poset $\cF(P)$. In fact,
as Bayer and Billera \cite{BB} showed,
they depend only on the flag numbers of $P$: given a sequence of integers 
$I = (i_1, \dots, i_n)$ with $0 \le i_1 < i_2 < \dots < i_n \le d$,
an $I$-flag is an $n$-tuple
$F_1 < F_2 < \dots < F_n$ of faces of $P$ with $\dim F_k = i_k$ for all $k$.
The $I$-th {\em flag number\/} $f^{}_I(P)$ is the number of $I$-flags.
Letting $P$ vary over all polytopes of a given dimension $d$,
the numbers $g_i(P)$ and $h_i(P)$ can be expressed as a $\Z$-linear
combination of the $f^{}_I(P)$.

Conjecturally all the $g_i(P)$ should be nonnegative for all $P$.
This is known to be true for $i = 1, 2$ \cite{K}.  For higher
values of $i$, it can be proved for rational polytopes using the
interpretation of $g_i(P)$ as an intersection homology Betti number 
of an associated toric variety.

\begin{prop} \label{gnneg} 
If $P$ is a rational polytope, then $g_i(P) \ge 0$ for all $i$.
\end{prop}

\section{Relative $g$-polynomials}

The following proposition defines a relative version of the classical
$g$-polynomials.
\begin{prop}\label{grdef} There is a unique family of polynomials $g(P,F)$ associated
to a polytope $P$ and a face $F$ of $P$, satisfying the following 
relation: for all $P, F$, we have
\begin{equation}\label{grdefeq}
 \sum_{F \le E \le P} g(E, F)g(P/E) = g(P).\end{equation}
\end{prop}
\begin{proof} The equation \eqref{grdefeq} can be used inductively to 
compute $g(P, F)$, since the left hand side gives $g(P, F)\cdot 1$ plus
terms involving $g(E, F)$ where $\dim E< \dim P$.  The induction starts
when $P = F$, which gives $g(F, F) = g(F)$.  Uniqueness is clear.
\end{proof}  
As an example, if $F$ is a facet of $P$, then
$g(P, F) = g(P) - g(F)$. 
Just as before we will denote the coefficient of
$q^i$ in $g(P, F)$ by $g_i(P,F)$.

We have the following notion of relative flag numbers.  Let $P$ be a 
$d$-polytope, and $F$ a face of dimension $e$.  Given a sequence of integers 
$I = (i_1, \dots, i_n)$ with $-1 \le i_1 < i_2 < \dots < i_n \le d$ and
a number $1 \le k \le n$ with $i_k \ge e$, define 
the relative flag number $f^{}_{I,k}(P, F)$ to be the number of $I$-flags
$(F_1, \dots,F_n)$ with $F \le F_k$.  Note that letting $k = n$ and
$i_n = d$ gives the ordinary flag numbers of $P$ as a special case.
Also note that the numbers $f^{}_{I, k}$ where $i_k = e$ give all 
products of the form $f^{}_J(F)f^{}_{J'}(P/F)$.

\begin{prop} Fixing $\dim P$ and $\dim F$, the relative $g$-number \\ 
$g_i(P, F)$ is a $\Z$-linear combination of the $f_{I, k}(P, F)$.
\end{prop}
\begin{proof} Use induction on $\dim P/F$.  If $P = F$, then we
have $g(P,P) = g(P)$ and the result is just the corresponding 
result for the ordinary flag numbers.  If $P \ne F$, the equation
\eqref{grdefeq} gives
\[ g(P, F) = g(P) - \sum_{e = \dim F}^{\dim P - 1}
\sum_{\substack{ \dim E = e \\ F\le E <  P}} g(E, F)g(P/E).\]
For every $e$ the coefficients of the 
inner summation on the right hand side are 
$\Z$-linear combinations of the $f_{I,k}(P, F)$, using the inductive hypothesis.
\end{proof}

The following theorem is the 
main result of this paper.  It will be a consequence of Theorem \ref{grint}.
\begin{thm} \label{grnneg} If $P$ is a rational polytope and $F$ is any face, then \\ $g_i(P, F) \ge 0$ for all $i$.
\end{thm}

\begin{cor}[Kalai's conjecture] If $P$ is a rational polytope and $F$ is
any face, then
\[g(P) \ge g(F)g(P/F),\]
where the inequality is taken coefficient by coefficient.
\end{cor}
\begin{proof} For any face $E$ of $P$ the polytope $P/E$ is rational,
 so we have $g(P) = g(F, F)g(P/F) \mathop{+}$ other nonnegative terms.
\end{proof}

\section{Some examples and formulas}
This section contains further combinatorial results on the relative
$g$-polynomials.  They are not used in the remainder of the paper.

First, we give an interpretation of $g_1(P, F)$ and $g_2(P, F)$
analogous to the ones Kalai gave for the usual $g_1$ and $g_2$ in
\cite{K}.  We begin by recalling those results from \cite{K}.

Give a finite set of points $V \subset \R^d$
define the space $\Aff(V)$ of affine dependencies of $V$ to be
\[\{\, a \in \R^V \mid \Sigma_{v\in V} a_v = 0,\, 
\Sigma_{v\in V} a_v\cdot v = 0\,\}.\] If $V_P$ is the set of vertices
of a polytope $P \subset \R^d$, then $\Aff(V_P)$ is a vector space
of dimension $g_1(P)$.

To describe $g_2(P)$ we need the notion of stress on a framework. 
A framework $\Phi = (V, E)$ 
is a finite collection $V$ of points in $\R^d$
together with a finite collection $E$ of straight line 
segments (edges) joining them.  Given a finite 
set $S$, we denote the standard basis 
elements of $\R^S$ by $1_s, s\in S$.  The space of stresses
$\cS(\Phi)$ is the kernel of the linear map
\[\alpha\colon \R^E \to \R^V \otimes \R^d,\] defined by
\[\alpha(1_e) = 1_{v_1}\otimes (v_1 - v_2) + 1_{v_2}\otimes (v_2 - v_1),\]
where $v_1$ and $v_2$ are the endpoints of the edge $e$.  A stress
can be described physically as  an assignment of a contracting
or expanding force to each edge, such that the total force resulting at 
each vertex is zero.

To a polytope $P$ we
can associate a framework $\Phi_P$ by taking as vertices the vertices of
$P$, and as edges the edges of $P$ together with enough extra edges to
triangulate all the $2$-faces of $P$.  Then $g_2(P)$ is the dimension
of $\cS(\Phi_P)$.

Given a polytope $P$ and a face $F$, define the closed union of faces
$N(P,F)$ to be the union of all facets of $P$ containing $F$.  
Note that $N(P,\emptyset) = \partial P$, and $N(P, P) = \emptyset$.
Let $V_N$ be the set of vertices of $P$ in $N(P, F)$, and define
a framework $\Phi_N$ by taking all edges and vertices of 
$\Phi_P$ contained in $N(P, F)$.

\begin{thm} We have
\[g_1(P, F) = \dim_\R \Aff(V_P)/ \Aff(V_N), \text{and}\]
\[g_2(P, F) = \dim_\R \cS(\Phi_P)/\cS(\Phi_N),\] using the 
obvious inclusions of $\Aff(V_N)$ in $\Aff(V_P)$ and $\cS(\Phi_N)$
in $\cS(\Phi_P)$.
\end{thm}
The proof for $g_1$ is an easy exercise; the proof for $g_2$ will
appear in a forthcoming paper \cite{BM}.

Next, we have a formula which shows that $g(P, F)$ can be decomposed
in the same way $g(P)$ was in Proposition \ref{grdef}.  Given two 
faces $E, F$ of a polytope $P$, let $E \vee F$ be the unique 
smallest face containing both $E$ and $F$.
\begin{prop} For any polytope $P$ and faces $F' \le F$ of $P$, we have
\[g(P, F) = \sum_{F'\le E} g(E, F')g(P/E, (E\vee F)/E).\]
\end{prop}
\begin{proof}
As usual, we show that this formula for $g(P, F)$ satisfies the defining
relation of Proposition \ref{grdef}.  Fix $F' \le F$, and define 
$\hat{g}(P, F)$ to be the above sum.  Then we have
\begin{eqnarray*}
\sum_{F\le D} \hat{g}(D,F)g(P/D) & = & 
\sum_{\substack{F'\le E \\ F\vee E\le D}}
g(P/D)g(E, F')g(D/E, (E\vee F)/E)\\
 &=& \sum_{F'\le E} g(E, F')g(P/E)\\
 & =& g(P).
\end{eqnarray*}
Since the computation of $g(P, F)$ from Proposition \ref{grdef} only
involves computation of $g(E, F)$ for other faces $E$ of $P$, this
proves that $\hat{g}(P,F) = g(P, F)$, as required.
\end{proof}

Finally, we can carry out the inversion implicit in Proposition
\ref{grdef} explicitly.  First we need the notion of polar polytopes.
Given a polytope $P \subset \R^d$, we can assume that
the origin lies in the interior of $P$ by moving $P$ by an affine motion.
The polar polytope $P^*$ is defined by
\[P^* = \{\,x \in (\R^*)^d \mid \langle x, y \rangle \le 1 \ \text{for all}\ 
y \in P\,\}.\]  The face poset $\cF(P^*)$ is canonically the opposite
poset to $\cF(P)$.  Define $\gb(P) = g(P^*)$.

\begin{prop} We have
\begin{equation}\label{eq1}
 g(P, F) = \sum_{F\le F' \le P} (-1)^{\dim P - \dim F'} 
g(F')\gb(P/F').\end{equation}
\end{prop}
\begin{proof} We use the following formula, due to Stanley \cite{St}:
For any polytope $P \ne \emptyset$, we have
\begin{equation}\label{eq2}
\sum_{\emptyset \le F \le P} (-1)^{\dim F} \gb(F)g(P/F) = 0.
\end{equation}
Now define $\hat{g}(P, F)$ to be the right hand side of \eqref{eq1}.  
We will show that
the defining property \eqref{grdefeq} of Proposition \ref{grdef} holds.

Pick a face $F$ of $P$.
We have, using \eqref{eq2},
\[
\sum_{F \le E \le P} g(E, F)g(P/E) = 
\sum_{F \le F' \le E \le P} (-1)^{\dim E - \dim F'}g(F')\gb(E/F')g(P/E) \]
\begin{eqnarray*} 
& = & \sum_{F\le F' \le P} g(F') \sum_{F' \le E \le P} (-1)^{\dim E - \dim F'}
\gb(E/F')g(P/E) \\
 & = & g(P),
\end{eqnarray*}
as required.
\end{proof}

\section*{Introduction to \S\S 4 - 6}
Sections 4-6 of this paper concern the topology of algebraic varieties.   
They may be read independently from the combinatorics of \S 1-3.   
 
The principal result, Theorem 16 of \S 5, is this:  Consider a subvariety  
$Y$ of a complex algebraic variety $X$.  Suppose that there is a blowup  
$p\colon\wt{X}\to X$ such that $p^{-1}Y$ has a neighborhood in 
$\wt{X}$ that is homeomorphic  
to a line bundle over $p^{-1}Y$.  Then the restriction of the intersection  
homology sheaf $\IC(X)$ of $X$ to $Y$ is a direct sum of shifted intersection  
homology sheaves.   
 
For our applications, we need a slight strengthening of this result.   
The neighborhood of $p^{-1}Y$ 
will only be a Seifert bundle $E\to B$, a generalization  
of a line bundle which allows fibers to be quotients by cyclic groups.   
These are treated in \S 4.  In \S 6 we apply the principal result to the 
inclusion of toric varieties $Y_F\subset X_P$, where $Y_F \simeq X_{P/F}$ 
is the closure of the torus orbit corresponding to $F$.  If $x$ is the 
unique torus-fixed point of $X$, then
$g_i(P,F)$ measures the number of copies of the intersection homology  
sheaf $\IC(\{x\})$  that appear with shift $2i$ in the restriction 
of the intersection homology sheaf of $X_P$ to $Y_F$.

\section{Seifert bundles}
In this section we investigate maps of algebraic varieties
 $E \to B$ which are nearly line bundles, but which allow fibers to
be quotients by cyclic groups.

\begin{defn} A {\em Seifert bundle \/} is an affine map 
$\pi\colon E \to B$ of algebraic
varieties, together with a section (which we will sometimes call
the zero section) $s\colon B \to E$, and an
algebraic $\C^*$-action on $E$, so that:
\begin{itemize}
\item giving $B$ the trivial $\C^*$-action, 
$\pi$ and $s$ are $\C^*$-equivariant,
\item each fiber $\pi^{-1}(b)$ is a curve whose normalization is isomorphic to 
the complex line, on which $\C^*$ acts by multiplication by a character
$x \mapsto x^{n_b}$, $n_b>0$.
\end{itemize}
\end{defn}

\begin{lemma} If $(E, B, \pi, s)$ is a Seifert bundle, and $n_b$ is a
constant on all of $B$, then $E$ is topologically a complex line
bundle over $B$.
\end{lemma}
\begin{proof} It is enough to show this locally, so assume $B$ is an affine
variety, and take $b\in B$.  Then $E$ is also affine, and the $\C^*$ action
on $E$ induces an action, and hence a nonnegative grading, 
on the coordinate ring
$A(E)$.  Take a polynomial $f$ which doesn't vanish on $\pi^{-1}(b)$;
we can assume it is homogeneous.  Shrinking $B$ if necessary, we can assume
$f$ doesn't vanish on any fiber of $\pi$.

Let $Y \subset E$ be the subvariety defined by the equation $(f = 1)$.
We will show that $\pi|_Y\colon Y \to B$ is proper.  The lemma
follows from this claim; if $\pi|_Y$ is proper, then the natural bijection
$Y/G \to B$, where $G$ is the group of $d$th roots of unity, $d = \deg f$,
 is proper and hence a homeomorphism.  Then, since the $G$-action is free, 
$Y \to B$ is a covering map.
For a small (topological) neighborhood 
$U$ of $b$ we can thus define a continuous section $\sigma\colon B\to Y$.
The map $(b, t) \mapsto t\cdot \sigma(b)$ gives the required local
trivialization of $E$.

To show the claim, take a compact set $K \subset B$.
Choosing a homogeneous system $(f_1, \dots, f_s)$ of 
generators for $A = A(E)$ over its zeroth graded piece $A_0$ defines
an embedding \[E \subset B \times \C^n\] as a closed subvariety, and the 
$\C^*$ action on $E$ is given by a linear $\C^*$ action on $\C^n$.
Let $r>0$ be the smallest character of $\C^*$ appearing in a 
diagonalization of this 
action; it is the smallest of the degrees of the $f_i$.
Then if $S^{2n-1}$ is the set of elements of norm one in $\C^n$,
the set \[\pi^{-1}(K) \cap (B \times S^{2n-1})\]
is compact, and so the values $|f|$ takes on it are bounded away from 
zero, say by $\delta$.  Thus $\pi^{-1}(K) \cap Y$ is a closed subset
of $K\times N_{1/\delta^r}$, where $N_a \subset \C^n$ 
is the closed ball of vectors of norm $\le a$, and so is compact.  

\end{proof}

\begin{cor} \label{Sbcor} Any Seifert bundle $E$ over $B$ maps to a 
(topological) line bundle
$E'$ over $B$ by a finite map.
\end{cor}
\begin{proof} We can take the least common multiple
\[n = \lcm_{b \in B} n_b\]
of the numbers $n_b$, since there are only finitely many distinct 
values of $n_b$.  Setting $E' = E/G$ where $G$ is the group of
$n$th roots of unity in $\C^*$ does the trick.
\end{proof}

Intuitively, the zero-section map $s$ for a Seifert bundle will be
a ``$\Q$-homology normally nonsingular inclusion".  We have the 
following generalization of a result of \cite{GM} about normally 
nonsingular inclusions:

\begin{prop}\label{Sbprop}
 Let $(E, B, \pi, s)$ be a Seifert bundle.  Then there is
an isomorphism
\[s^*\IC(E) \cong \IC(B).\]
\end{prop}

We need a small lemma first.

\begin{lemma}\label{fgqIC}
 Let $X$ be a pseudomanifold, acted on by a finite group $G$, and let
$Y$ be a $G$-invariant subspace.
Then there is an isomorphism
\[IH_*(X/G, Y/G; \Q) \cong IH_*(X, Y; \Q)^G\]
between the intersection homology of the pair 
$(X/G, Y/G)$ and the $G$-stable part
of the intersection homology of $(X, Y)$.
\end{lemma}
\begin{proof} Give $X$ a $G$-invariant triangulation.  Then the intersection
homology of $X$ can be expressed by means of simplicial chains of the
barycentric subdivision, see \cite[Appendix]{MV}.  Now the standard argument
in \cite[p. 120]{Br} can be applied.
\end{proof}

\begin{proof}[Proof of Proposition \ref{Sbprop}]
By Corollary \ref{Sbcor}, we can map $E$ to a line bundle $E'$ by a 
finite map.  Let $s'$, $\pi'$ denote the section and projection 
maps for $E'$.  Let $\mb{A} = s^*\IC(E)$, $\mb{A'} = (s')^*\IC(E')$.  Because
$E'$ is a line bundle over $B$, $\mb{A'}$ is isomorphic to $\IC(B)$.  

Let $U \subset B$ be a Zariski open subset where $n_b$ is constant.
Then $E|_U = \pi^{-1}(U)$ is a line bundle over $B$, so $F|_U$ is
a one-dimensional constant local system.  We will show that for any
point $p \in B$ there are isomorphisms
\[j_p^*\mb{A} \cong j_p^*\mb{A'},\ j_p^!\mb{A} \cong j_p^!\mb{A'}\]
between the stalks and costalks of $\mb{A}$ and $\mb{A'}$
(or more precisely isomorphisms between their cohomology groups),
where $j_p$ is the inclusion.  It follows that $F$ satisfies the 
perversity axioms defining the intersection homology sheaf from \cite{GM}.

To show the claim, note that since the $\C^*$ action retracts both $E$ and $E'$ 
onto $B$, we have isomorphisms
\[\mb{A} \cong R\pi_*\IC(E), \ \mb{A'} \cong R\pi'_*\IC(E').\]
So we can describe the stalks and costalks of $F$ and $F'$ as
follows.  Let $N$ be a small neighborhood of $p$ in $B$, and let $L$ be
its boundary.  Then we have
\begin{eqnarray*}
\HH^i j_p^*\mb{A} & = & IH_{n-i}(\pi^{-1}(N), \pi^{-1}(L); \Q),\\
\HH^i j^!_p\mb{A} & = & IH_{n-i}(\pi^{-1}(N); \Q),
\end{eqnarray*}
where $n$ is the real dimension of $B$, and similarly for $\mb{A'}$.  
The claim now follows from Lemma \ref{fgqIC},
using the fact that the finite group $G$ is contained in $\C^*$ and
hence acts trivially on the intersection homology groups above.
\end{proof}

\section{Seifert resolutions}

\begin{defn} A {\em Seifert resolution\/} of an inclusion $Y \subset X$
of irreducible algebraic varieties is a variety $\wt{X}$ together with 
a proper, surjective map
$p\colon \wt{X}\to X$, so that, if $\wt{Y} = p^{-1}(Y)$, then 
$p$ induces an isomorphism of $\wt{X} \setminus \wt{Y}$ with $X \setminus Y$,
and the inclusion $\wt{Y} \subset \wt{X}$ is the zero section of
a Seifert bundle.
\end{defn}

Now suppose $X$ is a connected normal algebraic variety with a nontrivial
 algebraic $\C^*$ action.  
Let $Y$ be an irreducible subvariety contained in the fixed point
set of $X$.  Let $U = \{\, x\in X \mid \overline{\C^*\cdot x} \cap Y 
\ne \emptyset\,\}$.  We say that $Y$ is an {\em attractor \/} for the
$\C^*$ action if for all points $x\in U$ the limit
 $\lim_{t\to 0} t\cdot x$ exists and lies in $Y$ and the only points
$x \in U$ for which $\lim_{t \to \infty} t\cdot x$ lies in $Y$ are
already in $Y$.

\begin{thm} \label{srex} If $U$ is an open neighborhood of
$Y$ and $Y$ is an attractor, then the pair $(Y, X)$ has a Seifert resolution.
\end{thm}
\begin{proof} 
We will show that $(Y, U)$ has a Seifert resolution; this will be enough.
By \cite{Su}, every point $y \in Y$ has a $\C^*$-invariant
affine neighborhood $U_y \subset U$.  Let $A_y$ be its coordinate 
ring.  The $\C^*$ action induces a grading on $A_y$ which is nonnegative 
because $Y$ is an attractor.  Further, if $R_y$ is the coordinate 
ring of $Y_y = U_y \cap Y$, the natural quotient map $A_y \to R_y$
identifies $R_y$ with the zeroth graded piece of $A_y$.  Thus there
is a projection map $\rho_y\colon U_y \to Y_y$; these glue to give
$\rho\colon U \to Y$.

Furthermore, the varieties and maps 
$\mathop{\rm Proj}(A_y) \to Y_y$ glue to give a variety $\wt{Y}$ and
a proper map $q\colon \wt{Y}\to Y$ (in other words, we let
$\wt{Y} = (U \setminus Y)/\C^*$). We also have a map 
$k\colon U\setminus Y \to \wt{Y}$ satisfying 
$q\circ k = \rho|_{U\setminus Y}$.

Define a morphism $U\setminus Y \to U\times_Y \wt{Y}$ by sending
$x$ to $(x, k(x))$.  Let $\wt{U}$ be the closure of the image of this
map, and let $p\colon \wt{U}\to U$ and $\pi\colon \wt{U}\to \wt{Y}$
be the restrictions of the projections of $U\times_Y \wt{Y}$ on the 
first and second factor, respectively.  $\wt{U}$ will be the required
Seifert resolution.

Note that $p^{-1}(Y) = Y \times_Y \wt{Y} \cong \wt{Y}$.  The
map $p$ is proper, because the projection \[U\times_Y \wt{Y} \to U \cong 
U \times_Y Y\] is proper.  It is now easy to check that $\wt{U}$ is a 
Seifert bundle over $\wt{Y}$.
\end{proof}

\begin{defn} Call an object $\mb{A}$ in $D^b(X)$ {\em pure}\/
if it is a direct sum of shifted intersection homology sheaves
\begin{equation}
\label{pureeq}\bigoplus_\alpha \IC(Z_\alpha; \cL_\alpha)[n_\alpha],
\end{equation}
where each $Z_\alpha$ is an irreducible subvariety of $X$, 
$\cL_\alpha$ is a simple
 local system on a Zariski open subset $U_\alpha$ of the smooth locus of 
$Z_\alpha$, and $n_\alpha$ is an integer.
\end{defn} 

\begin{lemma}\label{KSps}If $\mb{A, B}$ are objects in $D^b(X)$ and 
$\mb{A} \oplus \mb{B}$ is 
pure, then so is $\mb{A}$.
\end{lemma}
\begin{proof} Denote $\mb{A} \oplus \mb{B}$ by $\mb{C}$.  
Since $\mb{C}$ is pure, it is isomorphic to the direct sum 
\[\bigoplus_{i \in \Z} {}^p\!H^i(\mb{C})[-i]\]
of its perverse homology sheaves.  Each 
${}^p\!H^i(\mb{C}) = {}^p\!H^i(\mb{A}) \oplus {}^p\!H^i(\mb{B})$ 
is a pure perverse
sheaf, and since the category of perverse sheaves is abelian,
${}^p\!H^i(\mb{A})$ is pure.  Then the composition
\[\bigoplus {}^p\!H^i(\mb{A})[-i] \to \bigoplus {}^p\!H^i(\mb{C})[-i] 
\cong \mb{C} \to \mb{A}\]
induces an isomorphism on all the perverse homology sheaves, and
hence is an isomorphism (see \cite{BBD}).
\end{proof}

Also note that the decomposition \eqref{pureeq} of a pure object $\mb{A}$
is essentially unique: any other such decomposition will be the same
up to a reordering of the terms and replacing the local system $\cL_\alpha$
by another local system $\cL'_\alpha$ on $U'_\alpha$, so that 
$\cL_\alpha$ and $\cL'_\alpha$ agree on $U_\alpha \cap U'_\alpha$.

\begin{thm} \label{Srmain}
If a pair of varieties $(Y, X)$ has a Seifert resolution,
then the pullback $j^*\IC(X)$ of the intersection homology sheaf
by the inclusion is a pure object in $D^b(Y)$.
\end{thm}
\begin{proof} Consider the fiber square 
\[\begin{CD}
 \wt{Y} @>{\tilde{\jmath}}>> \wt{X}\\
 @VVqV                        @VVpV\\
 Y      @>j>>                 X\\
\end{CD}\]
where $j, \tilde{\jmath}$ are the inclusions, and $q = p|_{\wt{Y}}$.  
Because $p$ and $q$ are proper we have
\[Rq_*\tilde{\jmath}^*\IC(\wt{X}) \cong j^*Rp_*\IC(\wt{X}).\] 
The left hand side is $Rq_*\IC(\wt{Y})$ by Proposition \ref{Sbprop},
which is pure by the decomposition theorem of \cite{BBD}.
The decomposition theorem also implies that $\mb{A}=Rp_*\IC(\wt{X})$
is pure, and because $\wt{X} \to X$ is an isomorphism on a 
Zariski dense subset, the intersection homology sheaf of $X$
must occur in $\mb{A}$ with zero shift.
Thus the right hand side becomes
\[j^*(\IC(X)) \oplus j^*\mb{A}',\]
where $\mb{A}'$ is pure.  The result now follows from Lemma \ref{KSps}.
\end{proof}

\section{Toric varieties}
We will only sketch the properties of toric varieties that we will
need.  For a more complete presentation, see \cite{F}. 
Throughout this section let $P$ be a $d$-dimensional rational 
polytope in $\R^d$.

  Define a 
toric variety $X_P$ as follows.  Embed $\R^d$ into $\R^{d+1}$ by
\[(x_1, \dots, x_d) \mapsto (x_1, \dots, x_d, 1),\] and let $\sigma = 
\sigma^{}_P$ be the cone over the image of $P$ with apex at the origin
in $\R^{d+1}$.  It is a rational polyhedral cone with respect to the
standard lattice $N = \Z^{d+1}$.
More generally, if $F$ is a face of $P$, let
$\sigma^{}_F$ be the cone over the image of $F$; set 
$\sigma^{}_\emptyset = \{0\}$.

Then define $X = X_P$ to be the affine toric variety $X_\sigma$ 
corresponding to $\sigma$.  It is the variety 
$\mathop{\rm Spec} \C[M \cap \sigma^\vee]$,
where 
\[\sigma^\vee = \{ x \in (\R^{d+1})^* \mid \langle x, y \rangle \ge 0
\quad\text{for all}\quad y\in \sigma\,\}\] 
is the dual cone to $\sigma$, $M$ is the 
dual lattice to $N$, and $\C[M \cap \sigma^\vee]$ if the semigroup algebra 
of $M \cap \sigma^\vee$.  It is a $(d+1)$-dimensional
normal affine algebraic variety, on which the
torus $T = \Hom(M, \C^*)$ acts.  Let $f_v\colon X_P \to \C$ be the
regular function corresponding to the point $v \in M \cap \sigma^\vee$.

The orbits of the action of $T$ on $X$ are parametrized by the faces of $P$.
Let $F$ be any face of $P$, including the empty face, and let
\[\sigma^\bot_F = \{\, x\in \sigma^\vee \mid \langle x, y\rangle = 0 
\quad\text{for all}\quad y\in \sigma^{}_F\,\}\]
be the face of $\sigma^\vee$ dual to $\sigma^{}_F$.
Then the variety
\[O_F := \{\,x\in X \mid f_v(x) \ne 0 \iff v\in M\cap \sigma^\bot_F\,\}\]
is an orbit, isomorphic to the torus $(\C^*)^{d - e}$, where $e = \dim F$.
Furthermore, all $T$-orbits arise this way.  Thus $X_P$ has a unique 
$T$-fixed point $\{p\} = O_P$.

Given a face $F$, the union
\[ U_F = \bigcup_{E\le F} O_{E}\] is a $T$-invariant open neighborhood of $O_F$.
There is a non-canonical isomorphism $U_F \cong O_F \times X_F$ where
$X_F$ is the affine toric variety corresponding to $F$, considered as
a polytope in the affine space spanned by $F$, with the lattice given
by its intersection with $N$.  If $O'_E$ denotes
the orbit of $X_F$ corresponding to a face $E \le F$, then $O_E$ sits in 
$U_F \cong O_F \times X_F$ as $O_F \times O'_E$.

The union \[Y_F = \bigcup_{F\le E} O_{E}\] is the closure $\overline{O_F}$.
It is isomorphic to the affine toric variety $X_{P/F}$.  More precisely,
$Y_F$ is the affine toric variety corresponding to the cone 
$\tau = \sigma/\sigma_F$, the image of $\sigma$ projected into 
$\R^{d+1}/\mathop{\rm span}\sigma^{}_F$, with the lattice given by the 
image of $N$.  It is an easy exercise to 
show that $\tau$ is the cone over a polytope of the type $P/F$.

The connection between toric varieties and $g$-numbers of polytopes
is given by the following result.  Proofs appear in \cite{DL,Fie}.
\begin{prop} \label{BKM}The local intersection homology groups of $X_P$ are
described as follows.  Let $x$ be a point in $O_F$, and let 
$j_x$ be the inclusion.  Then
\[\dim \HH^{2i}j^*_x\IC(X_P) = g_i(F),\]
and $\HH^kj_x^*\IC(X_P)$ vanishes for odd $k$.
\end{prop}

Now fix a face $F$ of $P$.

\begin{lemma} There exists a $\C^*$ action
coming from a one-dimensional subtorus of $T$ so that
the fixed point set is $Y_F = \overline{O_F}$ and for any $x \in X_P$,
\[\lim_{t\to 0} t\cdot x \in Y_F.\]
\end{lemma}
\begin{proof} Let $a \in N\cap \sigma$ be a lattice point in the
relative interior of $\sigma^{}_F$.  This defines a $\C^*$ action on
$X_P$ by letting, for all $t\in \C^*$, $x\in X_P$, and 
$v \in M\cap \sigma^\vee$,
\[f_v(t \cdot x) = t^{\langle a, v\rangle}f_v(x).\]
The required property of this action is clear.
\end{proof}

Thus we can apply Theorem \ref{srex} to obtain a Seifert
resolution $\wt{Y}$ of the pair $(Y_F, X_P)$.  Although we will not need
this, a description of $\wt{Y}$ is quite interesting.  Let $\Delta(a)$
be the fan obtained by coning off all the faces of $\sigma$ to the
one-dimensional cone $\tau$ containing $a$.
Then $\wt{Y}$ is the toric variety $X_{\Delta(a)}$.

So by Theorem \ref{Srmain}, if $j\colon Y_F \to X_P$ is the inclusion,
the pullback $\mb{A} = j^*\IC(X_P)$ is a direct sum
\[\bigoplus_\alpha \IC(Z_\alpha, \cL_\alpha)[n_\alpha]\]
of shifted simple intersection homology sheaves.

\begin{lemma} All the terms in this decomposition are of the form
\[\IC(Y_{E}, \Q_{O_{E}})[n]\]
where $E \ge F$ and $-n$ is a nonnegative even integer.
\end{lemma}
\begin{proof} Since the sheaf $\IC(X_P)$ is invariant under the action
of $T$, so is the pullback $\mb{A}$; 
it follows that all the varieties $Z_\alpha$
are $T$-invariant.  Second, the isomorphism $U_F \cong O_F \times X_F$
implies that the homology sheaves of $\IC(X_P)$, and hence of $\mb{A}$,
are constant on each orbit; thus no nonconstant local systems can occur.
Finally, the assertion about the shifts follows from Proposition \ref{BKM}.
\end{proof}

Thus we can write
\begin{equation}
\label{Asplits}
\mb{A} = \bigoplus_{E\ge F} \bigoplus_{i\ge 0}\IC(Y_{E}; \Q)[-2i] \otimes 
V^i_{E},
\end{equation}
for some finite dimensional $\Q$-vector spaces $V^i_{E}$.  

Now we come to the main result, which gives an interpretation of the
combinatorially defined polynomials $g(P, F)$ for rational polytopes 
which implies nonnegativity, and hence Theorem \ref{grnneg}.  
Let $\{p\} = O_P$ be the unique 
$T$-fixed point of $X_P$.

\begin{thm} \label{grint}
The relative $g$-number $g_i(P,F)$ is given by
\[g_i(P, F) = \dim^{}_\Q V^i_P.\]
\end{thm}
\begin{proof} Taking this for the moment as a definition of $g(P,F)$, 
we will show that the defining relation of Proposition
\ref{grdef} holds.  First we need to interpret the vector spaces
$V^i_E$ for $F\le E \ne P$.  There is a commutative diagram of inclusions
\[\begin{CD}  Y'_F @>j'>> X_E \\
             @Vk'VV      @VVkV \\
              Y_F  @>j>> X_P
\end{CD}\]
where $j'$ be the inclusion of $Y'_F = \overline{O'_F}$ in $X_E$, 
$k$ is the inclusion of $X_E$ in $U_E \cong O_E \times X_E$ as
$\{x\} \times X_E$, and $k'$ is the restriction of $k$.  

Then $k$ is a normally nonsingular 
inclusion, so we have
\[(j')^*k^*\IC(X_P) = (j')^*\IC(X_{E}) =\]
\[\bigoplus_{F\le F'\le E} \bigoplus_{i\ge 0}\IC(Y'_{F'}; \Q)[-2i] \otimes 
W^i_{F'}\]
for some vector spaces $W^i_{F'}$. 
On the other hand, since $k'$ is a normally nonsingular 
inclusion, it is also equal to 
\[(k')^*\mb{A} = \bigoplus_{F\le F'\le E} 
\bigoplus_{i\ge 0}\IC(Y'_{F'}; \Q)[-2i] \otimes V^i_{F'}.\]  Comparing terms,
we see that $W^i_{F'} \cong V^i_{F'}$, so we have
\[\dim^{}_\Q V^i_E = g_i(E, F).\]

The theorem now follows; the defining relation of Proposition \ref{grdef}
expresses two different ways of writing the dimensions of the 
stalk intersection homology 
groups of $X_P$ at the fixed point $p$.  One the one hand,
they are given by the coefficients of $g(P)$, by Proposition \ref{BKM}.
On the other hand they are given by \[\sum_{F \le E \le P} g(E, F)g(P/E),\]
using \eqref{Asplits}.
\end{proof}

\end{document}